\documentclass{amsart}

\usepackage[T2A,OT1]{fontenc}
\usepackage[cp1251,utf8]{inputenc}
\usepackage[english,russian,ukrainian]{babel}
\usepackage{amsfonts,amsmath,amssymb,amsthm,enumerate}

\newcommand{\1}{1\!\!\,{\rm I}}

\begin{document}
\selectlanguage{english}

\title[Evolution of moments of isotropic Brownian stochastic flows]{Evolution of moments of isotropic Brownian stochastic flows}
\author{V.~V.~Fomichov}
\date{}

\begin{abstract}
In this paper we consider the asymptotic behaviour of all moments of the interparticle distance and of all mixed moments of an isotropic Brownian stochastic flow which serves as a smooth approximation of the Arratia flow.
\end{abstract}

\maketitle

\section{Introduction}

Studying the dynamics of the moments of different physical and chemical parameters is usually a crucial element of the stochastic analysis of such fields as chemical kinetics, subsurface solute transport, neural networks modelling, etc. For instance, in \cite{GrahamMcLaughlin1}, \cite{GrahamMcLaughlin2} some PDEs are obtained for the first and second unconditional and conditional moments of solute concentration (see also \cite[Chapter 4]{Dean}). The derivation of moment equations for $n$-dimensional Markov processes for any $n\in\mathbb{N}$ using derivative, or jump, moments is presented in \cite{SotiropoulosKaznessis}. Furthermore, relations between the moments of a given distribution may provide some information about it when its precise form is not known, which is often the case in real applications of stochastic analysis. For instance, the moments of centered jointly Gaussian random variables must satisfy the relations implied by the well-known Wick formula, which allows to establish that the distribution of a given random vector is not Gaussian and to determine how much it deviates from the Gaussian distribution. For a more detailed discussion of the moment approach in applications see \cite{SotiropoulosKaznessis} and the references in it.

In this paper we are concerned with the asymptotic behaviour of the moments of a smooth approximation of the Arratia flow. To be more precise, we consider the following stochastic integral equation:
$$
x(u,t)=u+\int_0^t\int_{\mathbb{R}}\varphi(x(u,s)-q)W(dq,ds),\quad t\geq 0,\quad u\in\mathbb{R},\eqno(1.1)
$$
where $W$ is a Wiener sheet on $\mathbb{R}\times\mathbb{R}_+$ and the function $\varphi\in C^\infty(\mathbb{R})$ has a compact support (we denote it as $\varphi\in C_0^\infty(\mathbb{R})$), is non-negative and such that $\varphi(q)=\varphi(-q)$, $q\in\mathbb{R}$, and $\|\varphi\|_{L_2(\mathbb{R})}^2:= \int_{\mathbb{R}}\varphi^2(q)dq=1$ (on integration with respect to a Wiener sheet see~\cite{Walsh}, \cite{Kotelenez}, \cite{Dorogovtsev2007}). Under these conditions on the function $\varphi$ this equation has a unique solution for every $u\in\mathbb{R}$, and the corresponding family of random mappings ${x(\cdot,t):\mathbb{R}\mapsto\mathbb{R}}$, $t\geq 0$, defines an isotropic Brownian stochastic flow of $C^\infty$-diffeomorphisms (see \cite{Kunita}). Moreover, if for every $\varepsilon>0$ we set
$$
\varphi_\varepsilon(q):=\dfrac{1}{\sqrt{\varepsilon}} \varphi\left( \dfrac{q}{\varepsilon}\right),\quad q\in\mathbb{R},
$$
and denote by $x_{\varepsilon}(u,t)$ the solution of the stochastic integral equation
$$
x_\varepsilon(u,t)=u+\int_0^t\int_{\mathbb{R}} \varphi_\varepsilon(x_\varepsilon(u,s)-q)W(dq,ds),\quad t\geq 0,\quad u\in\mathbb{R},
$$
then for all $n\in\mathbb{N}$ the $n$-point motion of the corresponding stochastic flow converges weakly in $C([0;1], \mathbb{R}^n)$ as $\varepsilon\rightarrow 0$ to the $n$-point motion of the Arratia flow (see \cite{Dorogovtsev2004}).

Here we consider the asymptotic behaviour of all moments of the interparticle distance and of all mixed moments of the stochastic flow generated by the solutions of equation $(1.1)$. Although we restrict ourselves to this special choice of stochastic flow, our main results can be easily transferred to a more general case. More precisely, Theorem 2.6 holds for any one-dimensional isotropic Brownian stochastic flow provided its covariance function $b$ satisfies the following conditions:
$$
b(z)=1\Longleftrightarrow z=0,\eqno(1.2)
$$
$$
\lim_{z\rightarrow 0}\dfrac{1-b(z)}{z^2}=\beta>0,\eqno(1.3)
$$
$$
\forall n\in\mathbb{N}:\quad \lim_{z\rightarrow \infty}z^nb(z)=0.\eqno(1.4)
$$
Note that, due to the non-negative definiteness of $b$, these conditions are sufficient for conditions $(2.7)$ and $(2.8)$ to be fulfilled. Moreover, they imply that
$$
\forall \delta>0:\quad \int_0^\delta\dfrac{zdz}{1-b(z)}=+\infty,
$$
and so, using Feller's criterion of accessibility, we conclude that the particles of the stochastic flow do not coalesce with probability one (see \cite{Matsumoto}). On the other hand, no additional conditions apart from the covariance function being continuous and equal to one at the point zero are needed for Theorem 3.3 to hold.

\section{Asymptotics of the interparticle distance}

As for any one-dimensional isotropic Brownian stochastic flow (see \cite{LeJan}), the distance between any two particles of the flow converges to zero almost surely:
$$
\forall u,v\in\mathbb{R}:\ \lim_{t\rightarrow +\infty}\left(x(u,t)-x(v,t)\right)=0\quad \text{a. s.}\eqno(2.1)
$$
In our case the behaviour of the stochastic process
$$
\xi_t\equiv\xi_t(u,v):=x(u,t)-x(v,t),\quad t\geq 0,
$$
can be described in more detail. We prove that for any $u,v\in\mathbb{R}$, $u>v$, the stochastic process $\{-\ln\xi_t\}_{t\geq 0}$ is almost surely of order $t$ as $t\rightarrow+\infty$ and that for any $u,v\in\mathbb{R}$, $u\neq v$, and $n\in\mathbb{N}$ the moment ${\bf E}\xi^n_t$ is of order $t^{\frac{n-1}{2}}$ as $t\rightarrow+\infty$.

However, before stating the exact results, for the sake of completeness, we formulate an elementary proposition that will be often used in this paper.

\vskip10pt
\noindent {\bf Proposition 2.1.} {\it Let $f$ be a locally (Riemann) integrable function on $[0;+\infty)$ such that
$$
\lim_{t\rightarrow +\infty}f(t)=c\in\mathbb{R}\cup\{-\infty,+\infty\}.
$$
Then
$$
\lim_{t\rightarrow +\infty}\dfrac{1}{t}\int_0^tf(s)ds=c.
$$
}

\vskip10pt
\noindent {\bf Proof.} First assume that $c\in\mathbb{R}$. Then for arbitrary $\varepsilon>0$ there exist $t_1>0$ and $t_2>0$ such that
$$
\forall t>t_1:\quad \left|f(t)-c\right|<\dfrac{\varepsilon}{2},
$$
$$
\forall t>t_2:\quad \left|\dfrac{1}{t}\int_0^{t_1}(f(s)-c)ds\right|< \dfrac{\varepsilon}{2}.
$$
Therefore, for any $t>\max\{t_1,t_2\}$ we have
$$
\left|\dfrac{1}{t}\int_0^tf(s)ds-c\right|\leq \left|\dfrac{1}{t}\int_0^{t_1}(f(s)-c)ds\right|+ \left|\dfrac{1}{t}\int_{t_1}^t(f(s)-c)ds\right|\leq
$$
$$
\leq\left|\dfrac{1}{t}\int_0^{t_1}(f(s)-c)ds\right|+ \dfrac{1}{t}\int_{t_1}^{t}\left|f(s)-c\right|ds< \dfrac{\varepsilon}{2}+\dfrac{t-t_1}{t}\cdot\dfrac{\varepsilon}{2}<\varepsilon,
$$
which proves the desired result for $c\in\mathbb{R}$.

Now assume that $c=+\infty$. Then for arbitrary $M>0$ there exist $t_1>0$, $t_2>0$, and $t_3>0$ such that
$$
\forall t>t_1:\quad f(t)>2M,
$$
$$
\forall t>t_2:\quad \left|\dfrac{1}{t}\int_0^{t_1}f(s)ds\right|< \dfrac{1}{2}M,
$$
$$
\forall t>t_3:\quad \dfrac{t-t_1}{t}>\dfrac{3}{4}.
$$
Therefore, for any $t>\max\{t_1,t_2,t_3\}$ we have
$$
\dfrac{1}{t}\int_0^tf(s)ds=\dfrac{1}{t}\int_0^{t_1}f(s)ds+ \dfrac{1}{t}\int_{t_1}^tf(s)ds>-\dfrac{1}{2}M+\dfrac{t-t_1}{t}\cdot 2M>M,
$$
which proves the desired result for $c=+\infty$.

The proof for $c=-\infty$ is similar. $\square$

\vskip10pt
Besides, to simplify our considerations we will need some notations. Set
$$
\Phi(z):=\int_{\mathbb{R}}\varphi(z+q)\varphi(q)dq,\quad z\in\mathbb{R}.
$$
It is easy to see that the function $\Phi$ has the following properties:
$$
\Phi\in C_0^\infty(\mathbb{R}),
$$
$$
\Phi(z)=\Phi(-z),\quad z\in\mathbb{R},
$$
$$
\forall z\in\mathbb{R}:\quad 0\leq \Phi(z)\leq 1;\quad \Phi(z)=1 \Longleftrightarrow z=0.
$$
The first two are obvious, and to prove the third, note (see \cite[Theorem 188]{HardyLittlewoodPolya}) that in H\"{o}lder's inequality
$$
\Phi(z)=\int_\mathbb{R}\varphi(z+q)\varphi(q)dq\leq \sqrt{\int_\mathbb{R}\varphi^2(z+q)dq}\cdot\sqrt{\int_\mathbb{R}\varphi^2(q)dq}=1
$$
the sign of equality is possible if and only if the functions $\varphi(z+\cdot)$ and $\varphi(\cdot)$ are proportional. This is now equivalent to the identity
$$
\varphi(z+q)\equiv\varphi(q),\quad q\in\mathbb{R},
$$
and since the function $\varphi$ has a compact support, it is possible if and only if $z=0$.

Note that the functions $z\mapsto \Phi(z)$ and $z\mapsto z^n\Phi(z)$, $n\in\mathbb{N}$, are bounded, and so, combining $(2.1)$ and the dominated convergence theorem yields the relations
$$
{\bf E} \Phi(x(u,t)-x(v,t))\rightarrow 1,\quad t\rightarrow+\infty,\eqno(2.2)
$$
$$
{\bf E} \left[\left(x(u,t)-x(v,t)\right)^n\Phi(x(u,t)-x(v,t))\right]\rightarrow 0,\quad t\rightarrow+\infty,\quad n\in\mathbb{N}.\eqno(2.3)
$$

Also, set
$$
\sigma(z):=\sqrt{2(1-\Phi(z))},\quad z\in\mathbb{R}.
$$
It follows then from the properties of $\Phi$ that
$$
\sigma\in C^1((-\infty;0])\cap C^1([0;+\infty)),
$$
$$
\sigma(z)=\sigma(-z),\quad z\in\mathbb{R},
$$
$$
0\leq\sigma(z)\leq \sqrt{2},\quad z\in\mathbb{R};\quad \sigma(z)=0 \Longleftrightarrow z=0.
$$

With the help of these functions the joint quadratic variation of $\{x(u,t)\}_{t\geq 0}$ and $\{x(v,t)\}_{t\geq 0}$ and the quadratic variation of $\{\xi_t\}_{t\geq 0}$ can be written as
$$
\left\langle x(u,\cdot), x(v,\cdot)\right\rangle_t=\int_0^t\int_{\mathbb{R}}\varphi(x(u,s)-q)\varphi(x(v,s)-q)dqds=
$$
$$
=\int_0^t\Phi(x(u,s)-x(v,s))ds,\quad t\geq 0,\quad \text{a. s.}
$$
and
$$
\left\langle\xi\right\rangle_t=\left\langle x(u,\cdot)-x(v,\cdot)\right\rangle_t= \int_0^t\int_{\mathbb{R}}\left[\varphi(x(u,s)-q)-\varphi(x(v,s)-q)\right]^2dqds=
$$
$$
=\int_0^t\left[2-2\Phi(x(u,s)-x(v,s))\right]ds= \int_0^t\sigma^2(x(u,s)-x(v,s))ds,\quad t\geq 0,\quad \text{a. s.}
$$

Now we can prove the following theorem (note that the right-hand side of the equality below coincides with the Lyapunov exponent of the stochastic flow).

\vskip10pt
\noindent {\bf Theorem 2.2.} {\it With probability one for any $u,v\in\mathbb{R}$, $u>v$,
$$
\lim_{t\rightarrow +\infty}\dfrac{1}{t}\ln\left(x(u,t)-x(v,t)\right)=-\dfrac{1}{2}L',
$$
where $L':=\|\varphi'\|_{L_2(\mathbb{R})}^2>0$.
}

\vskip10pt
\noindent {\bf Proof.} The stochastic process $\{\xi_t=x(u,t)-x(v,t)\}_{t\geq 0}$ is strictly positive for all $t\geq 0$ with probability one and satisfies the stochastic integral equation
$$
\xi_t=(u-v)+\int_0^t\sigma(\xi_s)d\beta_s,\quad t\geq 0,
$$
where $\{\beta_t\}_{t\geq 0}$ is a Wiener process defined on the same probability space (see the proof of Theorem 3.4). Therefore, It\^{o}'s formula yields that
$$
\dfrac{1}{t}\ln\xi_t=\dfrac{1}{t}\ln(u-v)+ \dfrac{1}{t}\int_0^t\dfrac{\sigma(\xi_s)}{\xi_s}d\beta_s- \dfrac{1}{2t}\int_0^t\dfrac{\sigma^2(\xi_s)}{\xi_s^2}ds,\quad t>0.
$$
On the one hand, it follows from $(2.1)$ and
$$
\lim_{z\rightarrow 0+}\dfrac{\sigma^2(z)}{z^2}= \lim_{z\rightarrow 0+}\dfrac{2(1-\Phi(z))}{z^2}=-\lim_{z\rightarrow 0+}\Phi''(z)=-\Phi''(0)=L'\eqno(2.4)
$$
that
$$
\lim_{t\rightarrow +\infty}\dfrac{\sigma^2(\xi_t)}{\xi_t^2}=L'\quad \text{a. s.},
$$
which implies that
$$
\lim_{t\rightarrow +\infty}\dfrac{1}{t}\int_0^t\dfrac{\sigma^2(\xi_s)}{\xi_s^2}ds=L'\quad \text{a. s.} \eqno(2.5)
$$
On the other hand, the martingale
$$
m_t:=\int_0^t\dfrac{\sigma(\xi_s)}{\xi_s}d\beta_s,\quad t\geq 0,
$$
has the quadratic variation
$$
\left\langle m\right\rangle_t=\int_0^t\dfrac{\sigma^2(\xi_s)}{\xi_s^2}ds,\quad t\geq 0,\quad \text{a. s.},
$$
which, due to the boundedness of the function $(0;+\infty)\ni z\mapsto \sigma^2(z)/z^2\in (0;+\infty)$, can be estimated in the following way:
$$
\left\langle m\right\rangle_t\leq Kt,\quad t\geq 0,\quad \text{a. s.},
$$
where $K:=\sup_{z>0} (\sigma^2(z)/z^2)$. So,
$$
\left\langle \dfrac{1}{\sqrt{K}}m\right\rangle_t\leq t,\quad t\geq 0,\quad \text{a. s.},
$$
and from the representation
$$
\dfrac{1}{\sqrt{K}}m_t=\widetilde{\beta}_{<\frac{1}{\sqrt{K}}m>_t},\quad t\geq 0,\quad \text{a. s.},
$$
where $\{\widetilde{\beta}_t\}_{t\geq 0}$ is a standard Wiener process (maybe defined on an extended probability space), we obtain
$$
\dfrac{1}{\sqrt{K}}m_t=\widetilde{\beta}_{<\frac{1}{\sqrt{K}}m>_t}\leq\sup_{0\leq s\leq t}\widetilde{\beta}_{<\frac{1}{\sqrt{K}}m>_s}\leq\max_{0\leq s\leq t}\widetilde{\beta}_s,\quad t\geq 0,\quad \text{a. s.},
$$
and
$$
\dfrac{1}{\sqrt{K}}m_t=\widetilde{\beta}_{<\frac{1}{\sqrt{K}}m>_t}\geq\inf_{0\leq s\leq t}\widetilde{\beta}_{<\frac{1}{\sqrt{K}}m>_s}\geq\min_{0\leq s\leq t}\widetilde{\beta}_s,\quad t\geq 0,\quad \text{a. s.}
$$
However, the law of the iterated logarithm implies that
$$
\lim_{t\rightarrow +\infty} \dfrac{1}{t}\min_{0\leq s\leq t}\widetilde{\beta}_s= \lim_{t\rightarrow +\infty} \dfrac{1}{t}\max_{0\leq s\leq t}\widetilde{\beta}_s=0\quad \text{a. s.},
$$
and so,
$$
\lim_{t\rightarrow +\infty} \dfrac{1}{t}m_t=0\quad \text{a. s.} \eqno(2.6)
$$
Thus, from $(2.5)$ and $(2.6)$ we conclude that for any $u,v\in\mathbb{R}$, $u>v$,
$$
\lim_{t\rightarrow +\infty}\dfrac{1}{t}\ln\left(x(u,t)-x(v,t)\right)=-\dfrac{1}{2}L'\quad \text{a. s.},
$$
and the assertion of the theorem now follows from the strict monotonicity of $x(u,t)$ with respect to the spatial variable for all $t\geq 0$ with probability one. $\square$

\vskip10pt
\noindent {\bf Corollary 2.3.} {\it With probability one for any $u,v\in\mathbb{R}$, $u\neq v$,
$$
\lim_{t\rightarrow +\infty}\dfrac{1}{t}\ln\left(1-\Phi(x(u,t)-x(v,t))\right)=-L'.
$$
}

\vskip10pt
\noindent {\bf Proof.} This follows immediately from $(2.4)$, Theorem 2.2 and the symmetry of the function $\Phi$. $\square$

\vskip10pt
To prove the next theorem we will need a result concerning the evolution of a mass distribution in an isotropic Brownian stochastic flow.

\vskip10pt
\noindent {\bf Theorem 2.4 (\cite[Chapter 3, Theorem 2.20]{Zirbel}).} {\it Let $\{F_{s,t}\}_{0\leq s\leq t< +\infty}$ be the isotropic Brownian stochastic flow generated by a Brownian motion $U$ in $C(\mathbb{R},\mathbb{R})$, i. e. such that
$$
F_{s,t}(u)=u+\int_s^tU(F_{s,r}(u),dr),\quad 0\leq s\leq t<+\infty,\quad u\in\mathbb{R}.
$$
Suppose that its covariance function $b$ has the following properties:
$$
\lim_{z\rightarrow 0}\dfrac{b_0-b(z)}{z^2}=\beta>0\eqno(2.7)
$$
and
$$
\|b^*-(b_0-b)\|^2<b_0\cdot\|b_0-b\|,\eqno(2.8)
$$
where $b_0:=b(0)$, $b^*$ is the minimal concave majorant of the function $b_0-b$ on $\mathbb{R}_+$ and
$$
\|f\|:=\sup_{z\in [0;+\infty)}\left|f(z)\right|.
$$
Also, suppose that a probability measure $M_0$ on $\mathbb{R}$ satisfies the condition
$$
\exists \varepsilon>0:\quad \int_{\mathbb{R}}e^{\varepsilon \left|u\right|} M_0(du)<+\infty.\eqno(2.9)
$$
Then for
$$
D_t=\int_{\mathbb{R}}(u-C_t)^2M_t(du),\quad t\geq 0,
$$
where
$$
M_t=M_0\circ F_{0,t}^{-1},\quad t\geq 0,
$$
and
$$
C_t=\int_{\mathbb{R}}uM_t(du),\quad t\geq 0,
$$
there exist functions $l$ and $m$ such that
$$
l(t)\leq\dfrac{1}{\sqrt{t}} {\bf E} D_t\leq m(t),\quad t>0,
$$
and
$$
\lim_{t\rightarrow +\infty} l(t)=c\cdot l_{\infty},
$$
$$
\lim_{t\rightarrow +\infty} m(t)=c\cdot m_{\infty},
$$
where $l_{\infty}$ and $m_{\infty}$ are strictly positive constants depending only on the function $b$, and the constant $c$ is given by
$$
c=\int_{\mathbb{R}}\int_{\mathbb{R}}\left|u-v\right| M_0(du) M_0(dv).
$$
}

\vskip10pt
\noindent {\bf Remark 2.5.} The constants $l_{\infty}$ and $m_{\infty}$ are defined by the equalities
$$
l_{\infty}=\sqrt{\dfrac{8}{\pi}}\cdot\left(\sqrt{\dfrac{\|b_0-b\|}{2}}-\dfrac{\|b^*-(b_0-b)\|}{\gamma}\right),
$$
$$
m_{\infty}=\sqrt{\dfrac{8}{\pi}}\cdot\dfrac{\|b_0-b\|}{\gamma},
$$
where $\gamma<\sqrt{2b_0}$ can be taken arbitrarily close to $\sqrt{2b_0}$ (certainly, with a change of the functions $l$ and $m$).

\vskip10pt
\noindent {\bf Theorem 2.6.} {\it For any $u,v\in\mathbb{R}$ and $n\in\mathbb{N}\cup \{0\}$
$$
\lim_{t\rightarrow +\infty} \dfrac{1}{t^n} {\bf E} \left(x(u,t)-x(v,t)\right)^{2n+1}= 2^n\cdot (2n+1)!!\cdot\left(u-v\right)\eqno(2.10)
$$
and
$$
c_*\cdot 2^n\cdot (2n+2)!!\cdot\left|u-v\right|\leq
\varliminf_{t\rightarrow +\infty} \dfrac{1}{t^{(2n+1)/2}} {\bf E} \left(x(u,t)-x(v,t)\right)^{2n+2}\leq
$$
$$
\leq\varlimsup_{t\rightarrow +\infty} \dfrac{1}{t^{(2n+1)/2}} {\bf E} \left(x(u,t)-x(v,t)\right)^{2n+2}\leq c^*\cdot 2^n\cdot (2n+2)!!\cdot\left|u-v\right|\eqno(2.11)
$$
with the constants $c_*$ and $c^*$ given by
$$
c_*=\dfrac{2}{\sqrt{\pi}}(1-\|F-(1-\Phi)\|)>0,
$$
where $F$ is the minimal concave majorant of the function $1-\Phi$ on $\mathbb{R}_+$, and
$$
c^*=\dfrac{2}{\sqrt{\pi}}.
$$
}

\vskip10pt
\noindent {\bf Proof.} Set
$$
h_m(t):= {\bf E} \left(x(u,t)-x(v,t)\right)^m,\quad m\in\mathbb{N}.
$$
Then, using It\^{o}'s formula and Fubini's theorem, we get
$$
h_{m+2}(t)={\bf E} \left(x(u,t)-x(v,t)\right)^{m+2}=\left(u-v\right)^{m+2}+
$$
$$
+(m+2)(m+1)\int_0^t{\bf E} \left[\left(x(u,s)-x(v,s)\right)^m(1-\Phi(x(u,s)-x(v,s)))\right]ds=
$$
$$
=\left(u-v\right)^{m+2}+(m+2)(m+1)\int_0^th_m(s)ds-
$$
$$
-(m+2)(m+1)\int_0^t{\bf E} \left[\left(x(u,s)-x(v,s)\right)^m\Phi(x(u,s)-x(v,s))\right]ds,
$$
and so, $(2.3)$ implies that
$$
h_{m+2}(t)=(m+2)(m+1)\int_0^th_m(s)ds+\overline{o}(t),\quad t\rightarrow +\infty.\eqno(2.12)
$$
Obviously, we have
$$
h_1(t)\equiv u-v,
$$
i. e. for $n=0$ relation $(2.10)$ is true. Assume that it is true for $n=k\geq 0$. Then for $n=k+1$ using l'H\^{o}pital's rule we get
$$
\lim_{t\rightarrow +\infty}\dfrac{h_{2k+3}(t)}{t^{k+1}}=(2k+3)(2k+2)\cdot\lim_{t\rightarrow +\infty} \dfrac{\int_0^th_{2k+1}(s)ds} {t^{k+1}}=
$$
$$
=2\cdot(2k+3)\cdot\lim_{t\rightarrow +\infty} \dfrac{h_{2k+1}(t)}{t^k}=2^{k+1}\cdot(2k+3)!!\cdot (u-v),
$$
and the principle of mathematical induction implies the first part.

To prove the second part we will apply Theorem 2.4. To verify the conditions note that $(2.4)$ implies $(2.7)$, and $(2.8)$ takes the form
$$
\|F-(1-\Phi)\|<1.
$$
However, it is easy to see that now we can construct a function $\hat{F}$ of the form
$$
\hat{F}(z)=\min\{1;\alpha z\},\quad z\geq 0,
$$
where the coefficient $\alpha>0$ is sufficiently large, such that
$$
1-\Phi(z)\leq \hat{F}(z),\quad z\geq 0,
$$
and so,
$$
\|F-(1-\Phi)\|\leq \|\hat{F}-(1-\Phi)\|<1
$$
(the last inequality follows from $(2.4)$ and the continuity of $\Phi$). Finally, if we set
$$
M_0=\dfrac{1}{2}(\delta_u+\delta_v),
$$
where $\delta_a$ is the Dirac measure at a point $a\in\mathbb{R}$, then condition $(2.9)$ is also fulfilled.
Thus, by Theorem 2.4, we have
$$
\dfrac{1}{2}\left|u-v\right|\cdot\sqrt{\dfrac{8}{\pi}}\cdot\left(\dfrac{1}{\sqrt{2}}-\dfrac{\|F-(1-\Phi)\|}{\gamma}\right)\leq \varliminf_{t\rightarrow +\infty} \dfrac{1}{4\sqrt{t}} {\bf E} (x(u,t)-x(v,t))^2\leq
$$
$$
\leq\varlimsup_{t\rightarrow +\infty} \dfrac{1}{4\sqrt{t}} {\bf E} (x(u,t)-x(v,t))^2\leq\dfrac{1}{2}\left|u-v\right|\cdot\sqrt{\dfrac{8}{\pi}}\cdot\dfrac{1}{\gamma}.
$$
Recalling that $\gamma$ can be chosen arbitrarily close to $\sqrt{2}$, we get
$$
c_*\cdot 2\left|u-v\right|\leq \varliminf_{t\rightarrow +\infty} \dfrac{1}{\sqrt{t}} h_2(t)\leq \varlimsup_{t\rightarrow +\infty} \dfrac{1}{\sqrt{t}} h_2(t)\leq c^*\cdot 2\left|u-v\right|
$$
with the constants $c_*$ and $c^*$ defined above, i. e. for $n=0$ relation $(2.11)$ is true. Assuming that it is true for $n=k\geq 0$ we obtain that for arbitrary constants
$$
c_1<c_*\cdot 2^k\cdot(2k+2)!!\cdot\left|u-v\right|
$$
and
$$
c_2>c^*\cdot 2^k\cdot(2k+2)!!\cdot\left|u-v\right|
$$
there exists $t_0>0$ such that for any $t>t_0$ the inequalities
$$
c_1\leq \dfrac{h_{2k+2}(t)}{t^{(2k+1)/2}}\leq c_2
$$
hold. Hence, for any $t>t_0$ we have
$$
\int_0^th_{2k+2}(s)ds=\overline{o}(t)+\int_{t_0}^th_{2k+2}(s)ds\geq
$$
$$
\geq\overline{o}(t)+\dfrac{2c_1}{2k+3}\cdot (t^{(2k+3)/2}-t_0^{(2k+3)/2})=
$$
$$
=\overline{o}(t)+\dfrac{2c_1}{2k+3}\cdot t^{(2k+3)/2},\quad t\rightarrow +\infty,
$$
and
$$
\int_0^th_{2k+2}(s)ds=\overline{o}(t)+\int_{t_0}^th_{2k+2}(s)ds\leq
$$
$$
\leq\overline{o}(t)+\dfrac{2c_2}{2k+3}\cdot (t^{(2k+3)/2}-t_0^{(2k+3)/2})=
$$
$$
=\overline{o}(t)+\dfrac{2c_2}{2k+3}\cdot t^{(2k+3)/2},\quad t\rightarrow +\infty.
$$
So, using $(2.12)$ we get
$$
2c_1\cdot (2k+4)\leq\varliminf_{t\rightarrow +\infty} \dfrac{1}{t^{(2k+3)/2}} {\bf E} \left(x(u,t)-x(v,t)\right)^{2k+4}\leq
$$
$$
\leq\varlimsup_{t\rightarrow +\infty} \dfrac{1}{t^{(2k+3)/2}} {\bf E} \left(x(u,t)-x(v,t)\right)^{2k+4}\leq 2c_2\cdot (2k+4).
$$
Recalling that the constants $c_1$ and $c_2$ can be taken arbitrarily close to their bounds, we conclude that $(2.11)$ is also true for $n=k+1$. Applying the principle of mathematical induction completes the proof. $\square$

\section{Asymptotics of the mixed moments}

In this section we prove that for any $n\in\mathbb{N}$ and $u_1,\ldots,u_{2n}\in\mathbb{R}$ the mixed moment ${\bf E} \left[x(u_1,t)\ldots x(u_{2n},t)\right]$ is of order $t^n$ as $t\rightarrow+\infty$ and that for any $n\in\mathbb{N}$ and $u_1,\ldots,u_{2n-1}\in\mathbb{R}$ the mixed moment ${\bf E} \left[x(u_1,t)\ldots x(u_{2n-1},t)\right]$ is $\bar{o}(t^{n-\frac{1}{2}})$ as $t\rightarrow+\infty$.

We will need some additional notations. For $n\in\mathbb{N}$ denote by $C^n$ the space $C([0;1],\mathbb{R}^n)$ and define the norm
$$
\|\vec{f}\|_n:=\max_{1\leq k\leq n}\max_{0\leq t\leq 1} \left|f_k(t)\right|,
$$
where $\vec{f}=(f_1,\ldots,f_n)\in C^n$. Note that $C^n$ with the metric induced by this norm is a complete separable metric space.

Also, for any $u\in\mathbb{R}$ set
$$
\overline{x}(u,t):=x(u,t)-u,\quad t\geq 0,
$$
and for any $T>0$ and $u\in\mathbb{R}$ set
$$
\overline{x}_T(u,t):=\dfrac{1}{\sqrt{T}}\overline{x}(u,Tt)\equiv\dfrac{1}{\sqrt{T}}(x(u,Tt)-u),\quad 0\leq t\leq 1.
$$
Finally, $\vec{x}_T$ will stand for the random element $(\overline{x}_T(u_1,\cdot),\ldots,\overline{x}_T(u_n,\cdot))$ in the space $C^n$ for arbitrary $u_1,\ldots,u_n\in\mathbb{R}$.

\vskip10pt
\noindent {\bf Lemma 3.1.} {\it The following propositions are true:
\begin{flushleft}
$
({\it i})\ \forall u, v\in \mathbb{R}\quad\forall\varepsilon>0:\ \lim_{T\rightarrow+\infty}{\bf P}\{\max_{0\leq t\leq 1} \left|\overline{x}_{T}(u,t)-\overline{x}_{T}(v,t)\right|>\varepsilon\}=0;
$
$
({\it ii})\ {\bf P}\{\forall u, v\in \mathbb{R}:\ \lim_{T\rightarrow+\infty}\max_{0\leq t\leq 1} \left|\overline{x}_{T}(u,t)-\overline{x}_{T}(v,t)\right|=0\}=1;
$
$
({\it iii})\ \forall R>0:\ \lim_{T\rightarrow+\infty}\sup_{u,v\in[-R;R]}{{\bf E}\left(\max_{0\leq t\leq 1} \left|\overline{x}_{T}(u,t)-\overline{x}_{T}(v,t)\right|\right)^2}=0.
$
\end{flushleft}
}

\vskip10pt
\noindent {\bf Proof.} Clearly, $({\it i})$ is implied by $({\it ii})$, and $(ii)$ follows from the relations
$$
0\leq\max_{0\leq t\leq 1} \left|\overline{x}_{T}(u,t)-\overline{x}_{T}(v,t)\right|=
$$
$$
=\max_{0\leq t\leq 1} \left|\dfrac{1}{\sqrt{T}}\left(x(u,Tt)-x(v,Tt)\right)-\dfrac{1}{\sqrt{T}}\left(u-v\right)\right|\leq
$$
$$
\leq\dfrac{1}{\sqrt{T}}\max_{0\leq t\leq T} \left|x(u,t)-x(v,t)\right|+ \dfrac{1}{\sqrt{T}}\left|u-v\right|
$$
and $(2.1)$ together with the monotonicity of $x(u,t)$ with respect to the spatial variable for all $t\geq 0$ with probability one.

\noindent To prove $(iii)$ note that from Doob's inequality we get the estimate
$$
{\bf E}\left(\max_{0\leq t\leq 1} \left|\overline{x}_{T}(u,t)-\overline{x}_{T}(v,t)\right|\right)^2\leq
4{\bf E}\left(\overline{x}_{T}(u,1)-\overline{x}_{T}(v,1)\right)^2=
$$
$$
=4\left({\bf E}\overline{x}^2_{T}(u,1)+{\bf E}\overline{x}^2_{T}(v,1)-2{\bf E}\overline{x}_{T}(u,1)\overline{x}_{T}(v,1)\right)=
$$
$$
=8\left(1-\dfrac{1}{T}{\bf E}\overline{x}(u,T)\overline{x}(v,T)\right)= \dfrac{4}{T}\int_0^{T}{\bf E}\sigma^2(x(u,s)-x(v,s))ds.
$$
Using monotonicity again, we obtain that with probability one
$$
\sigma(x(u,t)-x(v,t))\leq\sigma^*(x(u,t)-x(v,t))\leq\sigma^*(x(R,t)-x(-R,t)),\quad t\geq 0,
$$
where
$$
\sigma^*(z):=\sup_{z'\in[-z;z]}\sigma(z'),\quad z\in\mathbb{R}_+,
$$
and so, we can write
$$
{\bf E}\left(\max_{0\leq t\leq 1} \left|\overline{x}_{T}(u,t)-\overline{x}_{T}(v,t)\right|\right)^2\leq \dfrac{4}{T}\int_0^{T}{\bf E}\sigma^{*2}(x(R,s)-x(-R,s))ds.\eqno(3.1)
$$
Combining $(3.1)$ and
$$
{\bf E}\sigma^{*2}(x(R,t)-x(-R,t))\rightarrow\sigma^{*2}(0)=0, \quad t\rightarrow+\infty,
$$
which is implied by the dominated convergence theorem, yields the required result. $\square$

\vskip10pt
\noindent {\bf Lemma 3.2.} {\it Let $\varkappa_T$ be the distribution of $\vec{x}_T$ in $C^n$, $\varkappa_w$ be the distribution of the random element $\vec{w}=(w(\cdot),\ldots,w(\cdot))$ in $C^n$, where $\{w(t)\}_{t\in [0;1]}$ is a standard Wiener process. Then $\varkappa_T$ converges weakly to $\varkappa_w$ as $T\rightarrow+\infty$.}

\vskip10pt
\noindent {\bf Proof.} Because of the scaling invariance of the Wiener process the marginal distributions of every measure $\varkappa_T$ coincide with the distribution of a standard Wiener process. Hence, the family of probability measures $\{\varkappa_T\}_{T>0}$ is weakly compact. Therefore, it is enough to show that for any sequence $\{T_k\}_{k=1}^\infty$ of positive real numbers, for which
$$
\lim_{k\rightarrow\infty}T_k=+\infty
$$
and the weak limit
$$
\lim_{k\rightarrow\infty}\varkappa_{T_k}=:\varkappa
$$
exists, the equality
$$
\varkappa=\varkappa_w
$$
holds. To do this, note that by Lemma 3.1 for any $\varepsilon>0$ and $i,j\in\{1,\ldots,n\}$ we have
$$
0\leq \int_{C^n} \1 \{\max_{0\leq t\leq 1}|f_i(t)-f_j(t)|>\varepsilon\}\varkappa(d\vec{f})\leq
$$
$$
\leq\varliminf_{k\rightarrow\infty}\int_{C^n} \1 \{\max_{0\leq t\leq 1}|f_i(t)-f_j(t)|>\varepsilon\}\varkappa_{T_k}(d\vec{f})=
$$
$$
=\varliminf_{k\rightarrow\infty}{\bf P}\{\max_{0\leq t\leq 1} \left|\overline{x}_{T_k}(u_i,t)-\overline{x}_{T_k}(u_j,t)\right|>\varepsilon\}=0,
$$
where $\1 \{A\}$ stands for the indicator function of a set $A$. This implies that
$$
\varkappa(\{\vec{f}\in C^n\ \vert\ f_1=\ldots=f_n\})=1.
$$
To finish the proof note that the marginal distributions of $\varkappa$ also coincide with the distribution of a standard Wiener process. $\square$

\vskip10pt
\noindent {\bf Theorem 3.3.} {\it The following propositions are true:
$$
\forall n\geq 1\quad\forall u_1, \dots, u_{2n-1}\in\mathbb{R}:\ \lim_{t\rightarrow +\infty}\dfrac{1}{t^{n-\frac{1}{2}}}{\bf E}\left[x(u_1,t)\ldots x(u_{2n-1},t)\right]=0,\eqno(3.2)
$$
$$
\forall n\geq 1\quad\forall u_1, \dots, u_{2n}\in\mathbb{R}:\ \lim_{t\rightarrow +\infty}\dfrac{1}{t^n}{\bf E}\left[x(u_1,t)\ldots x(u_{2n},t)\right]=(2n-1)!!.\eqno(3.3)
$$
}

\vskip10pt
\noindent {\bf Proof.} Note that for every $p>0$ we have
$$
\sup_{T>0}\int_{C^n}\|\vec{f}\|_n^p\varkappa_T(d\vec{f})= \sup_{T>0} {\bf E} \left[\max_{1\leq i\leq n}\max_{0\leq t\leq 1}\left|\overline{x}_T(u_i,t)\right|^p\right]\leq
$$
$$
\leq\sup_{T>0} \sum_{i=1}^{n} {\bf E} \left[\max_{0\leq t\leq 1}\left|\overline{x}_T(u_i,t)\right|^p\right]=n\cdot{\bf E} \left[\max_{0\leq t\leq 1}\left|w(t)\right|^p\right]<+\infty,
$$
where $\{w(t)\}_{t\in [0;1]}$ is a standard Wiener process. So, for any $s\geq0$
$$
\lim_{T\rightarrow +\infty} {\bf E} \left[\overline{x}_T(u_1,s)\ldots \overline{x}_T(u_{2n-1},s)\right]=
$$
$$
=\lim_{T\rightarrow +\infty} \int_{C^{2n-1}}\delta_s(f_1\cdot\ldots\cdot f_{2n-1})\varkappa_T(d\vec{f})=
$$
$$
=\int_{C^{2n-1}}\delta_s(f_1\cdot\ldots\cdot f_{2n-1})\varkappa_w(d\vec{f})= {\bf E}\left(w(s)\right)^{2n-1}=0
$$
and
$$
\lim_{T\rightarrow +\infty} {\bf E} \left[\overline{x}_T(u_1,s)\ldots \overline{x}_T(u_{2n},s)\right]=
$$
$$
=\lim_{T\rightarrow +\infty} \int_{C^{2n}}\delta_s(f_1\cdot\ldots\cdot f_{2n})\varkappa_T(d\vec{f})=
$$
$$
=\int_{C^{2n}}\delta_s(f_1\cdot\ldots\cdot f_{2n})\varkappa_w(d\vec{f})= {\bf E}\left(w(s)\right)^{2n}=(2n-1)!!\cdot s^n,
$$
where $\delta_s$ is the delta function at the point $s$. On the other hand, for any $s>0$
$$
\lim_{T\rightarrow +\infty} {\bf E} \left[\overline{x}_T(u_1,s)\ldots \overline{x}_T(u_{2n-1},s)\right]=
$$
$$
=\lim_{T\rightarrow +\infty} \dfrac{1}{T^{n-\frac{1}{2}}} {\bf E} \left[\overline{x}(u_1,Ts)\ldots \overline{x}(u_{2n-1},Ts)\right]=
$$
$$
=s^{n-\frac{1}{2}}\cdot\lim_{t\rightarrow +\infty} \dfrac{1}{t^{n-\frac{1}{2}}}{\bf E} \left[\overline{x}(u_1,t)\ldots \overline{x}(u_{2n-1},t)\right]
$$
and
$$
\lim_{T\rightarrow +\infty} {\bf E} \left[\overline{x}_T(u_1,s)\ldots \overline{x}_T(u_{2n},s)\right]=
$$
$$
=\lim_{T\rightarrow +\infty} \dfrac{1}{T^n}{\bf E} \left[\overline{x}(u_1,Ts)\ldots \overline{x}(u_{2n},Ts)\right]=
$$
$$
=s^n\cdot\lim_{t\rightarrow +\infty} \dfrac{1}{t^n}{\bf E} \left[\overline{x}(u_1,t)\ldots \overline{x}(u_{2n},t)\right].
$$
Thus, we obtain that
$$
\forall n\geq 1\quad\forall u_1, \dots, u_{2n-1}\in\mathbb{R}:\ \lim_{t\rightarrow +\infty}\dfrac{1}{t^{n-\frac{1}{2}}}{\bf E}\left[\overline{x}(u_1,t)\ldots \overline{x}(u_{2n-1},t)\right]=0\eqno(3.4)
$$
and
$$
\forall n\geq 1\quad\forall u_1, \dots, u_{2n}\in\mathbb{R}:\ \lim_{t\rightarrow +\infty}\dfrac{1}{t^n}{\bf E}\left[\overline{x}(u_1,t)\ldots \overline{x}(u_{2n},t)\right]=(2n-1)!!.\eqno(3.5)
$$

To prove $(3.2)$ and $(3.3)$ we will use the principle of mathematical induction. Note that for any $u,v\in\mathbb{R}$ we have
$$
\dfrac{1}{\sqrt{t}}{\bf E} x(u,t)=\dfrac{u}{\sqrt{t}}\rightarrow 0,\quad t\rightarrow +\infty,
$$
and
$$
\dfrac{1}{t}{\bf E} \left[x(u,t)x(v,t)\right]=\dfrac{1}{t}{\bf E} \left[\overline{x}(u,t)\overline{x}(v,t)\right]+\dfrac{uv}{t}=
$$
$$
=\dfrac{1}{t}\int_0^t{\bf E}\Phi(x(u,s)-x(v,s))ds+\dfrac{uv}{t}\rightarrow 1,\quad t\rightarrow +\infty,
$$
i. e. for $n=1$ both $(3.2)$ and $(3.3)$ are true. Assume that they are true for $n=k\geq 1$. Then for $n=k+1$ we have
$$
\overline{o}(1)=\dfrac{1}{t^{k+\frac{1}{2}}}{\bf E} \left[\overline{x}(u_1,t)\ldots \overline{x}(u_{2k+1},t)\right]=
$$
$$
=\dfrac{1}{t^{k+\frac{1}{2}}}{\bf E} \left[(x(u_1,t)-u_1)\ldots (x(u_{2k+1},t)-u_{2k+1})\right]=
$$
$$
=\dfrac{1}{t^{k+\frac{1}{2}}}\left({\bf E} \left[x(u_1,t)\ldots x(u_{2k+1},t)\right]+\overline{o}(t^{k+\frac{1}{2}})\right)=
$$
$$
=\dfrac{1}{t^{k+\frac{1}{2}}} {\bf E} \left[x(u_1,t)\ldots x(u_{2k+1},t)\right]+\overline{o}(1),\quad t\rightarrow +\infty,
$$
and so,
$$
\dfrac{1}{t^{k+1}}{\bf E} \left[\overline{x}(u_1,t)\ldots \overline{x}(u_{2k+2},t)\right]=
$$
$$
=\dfrac{1}{t^{k+1}}{\bf E} \left[(x(u_1,t)-u_1)\ldots (x(u_{2k+2},t)-u_{2k+2})\right]=
$$
$$
=\dfrac{1}{t^{k+1}}\left({\bf E} \left[x(u_1,t)\ldots x(u_{2k+2},t)\right]+\overline{o}(t^{k+\frac{1}{2}})\right)=
$$
$$
=\dfrac{1}{t^{k+1}} {\bf E} \left[x(u_1,t)\ldots x(u_{2k+2},t)\right]+\overline{o}(1),\quad t\rightarrow +\infty,
$$
which implies that $(3.2)$ and $(3.3)$ are also true for $n=k+1$. Applying the principle of mathematical induction yields the desired result. $\square$

\vskip10pt
Although Theorem 3.3 establishes the exact asymptotic behaviour of all even moments, the result concerning the odd moments is not the best possible. This is shown by Proposition 3.10. Its proof is based on the following theorem, which itself can be of interest.

\vskip10pt
\noindent {\bf Theorem 3.4.} {\it For any $u,v\in\mathbb{R}$ and $t\geq 0$ the equalities
$$
{\bf E} \left(\overline{x}(u,t)\vert\overline{x}(u,s)-\overline{x}(v,s), 0\leq s\leq t\right) = \dfrac{1}{2}\left(\overline{x}(u,t)-\overline{x}(v,t)\right),\eqno(3.6)
$$
$$
{\bf E} \left(\overline{x}(v,t)\vert\overline{x}(u,s)-\overline{x}(v,s), 0\leq s\leq t\right) = -\dfrac{1}{2}\left(\overline{x}(u,t)-\overline{x}(v,t)\right)\eqno(3.7)
$$
hold almost surely.
}

\vskip10pt
For the proof of this theorem we need the following two results.

\vskip10pt
\noindent {\bf Theorem 3.5 (\cite[Chapter 5, Theorem 5.12]{LiptserShiryaev}).} {\it Suppose that the martingale $m=\left(m_t,\mathcal{F}_t\right)_{t\in[0;T]}$ has continuous trajectories and its quadratic variation can be represented in the form
$$
\left\langle m\right\rangle_t=\int_0^ta_s^2ds,\quad 0\leq t\leq T,
$$
where the nonanticipative function $a_t=a(\omega,t)$ is such that $a(\omega,t)>0$ almost everywhere on $\Omega\times[0;T]$ with respect to the measure ${\bf P}\otimes \lambda$, where $\lambda$ is the one-dimensional Lebesgue measure. Then, on the initial probability space, there exists a standard Wiener process $\beta=\left(\beta_t,\mathcal{F}_t\right)_{t\in[0;T]}$ such that with probability one
$$
m_t=m_0+\int_0^ta_sd\beta_s,\quad 0\leq t\leq T.
$$
}

\vskip10pt
\noindent {\bf Remark 3.6.} The Wiener process $\{\beta_t\}_{t\in[0;T]}$ in Theorem 3.5 can be defined as
$$
\beta_t=\int_0^t\dfrac{dm_s}{a_s},\quad 0\leq t\leq T.
$$

\vskip10pt
\noindent {\bf Lemma 3.7.} {\it Let the stochastic process $\eta=(\eta_t,\mathcal{F}_t)_{t\in [0;T]}$ be a strong solution of the stochastic integral equation
$$
\eta_t=\eta_0+\int_0^tb(\eta_s)d\beta_s,\quad 0\leq t\leq T,
$$
where $\beta=(\beta_t,\mathcal{F}_t)_{t\in [0;T]}$ is a standard Wiener process and the function $b$ is such that
$$
\left|b(u)\right|\leq C\cdot(1+\left|u\right|),\quad u\in\mathbb{R},
$$
$$
{\bf P} \left\lbrace\int_0^Tb^2(\eta_t)dt<+\infty\right\rbrace =1,
$$
and
$$
b(\eta_t)>0
$$
almost everywhere on $\Omega\times[0;T]$ with respect to the measure ${\bf P}\otimes \lambda$, where $\lambda$ is the one-dimensional Lebesgue measure.
	
Then for any square-integrable martingale $m=(m_t,\mathcal{F}_t)_{t\in [0;T]}$ the stochastic process $\{m_t^{\eta}={\bf E} (m_t|\mathcal{F}_t^{\eta})\}_{t\in [0;T]}$, where $\mathcal{F}_t^{\eta}:=\sigma(\eta_s, 0\leq s\leq t)$, has the following properties:
\begin{flushleft}
1) $(m_t^{\eta},\mathcal{F}_t^{\eta})_{t\in [0;T]}$ is a square-integrable martingale;

2) $(m_t^{\eta},\mathcal{F}_t^{\eta})_{t\in [0;T]}$ has a continuous modification;

3) the continuous modification of $(m_t^{\eta},\mathcal{F}_t^{\eta})_{t\in [0;T]}$ permits {\bf P}-a. s. the representation
$$
m_t^{\eta}=m_0^{\eta}+\int_0^t{\bf E} \left(\dfrac{d}{ds}\left\langle m,\beta\right\rangle_s|\mathcal{F}_s^{\eta}\right)d\beta_s,\quad 0\leq t\leq T.
$$
\end{flushleft}
}

\vskip10pt
\noindent The {\bf proof} is similar to that of \cite[Chapter 5, Theorem 5.16]{LiptserShiryaev} and \cite[Chapter 8, Theorem 8.1]{LiptserShiryaev} and therefore omitted. $\square$

\vskip10pt
\noindent {\bf Remark 3.8.} Theorem 3.5 and Lemma 3.7 can be easily extended to the case of the infinite time interval $[0;+\infty)$ (in Lemma 3.7 square-integrability only on finite intervals should be used in both places). It is in this form that we will use them in the proof of Theorem 3.4.

\vskip10pt
\noindent {\bf Proof of Theorem 3.4.} If $u=v$, then equalities $(3.6)$ and $(3.7)$ are obvious. Therefore, we assume that $u\neq v$. Then the diffeomorphic property of the stochastic flow implies that with probability one
$$
\xi_t=x(u,t)-x(v,t)\neq 0,\quad t\geq 0,
$$
and so, with probability one
$$
\sigma(\xi_t)>0,\quad t\geq 0.
$$
Thus, by Theorem 3.5, we can write
$$
\xi_t=(u-v)+\int_0^t\sigma(\xi_s)d\beta_s,\quad t\geq 0,\quad \text{a. s.},
$$
where
$$
\beta_t=\int_0^t\dfrac{d\xi_s}{\sigma(\xi_s)},\quad t\geq 0.
$$
Then
$$
\left\langle\overline{x}(u,\cdot),\beta\right\rangle_t=\left\langle x(u,\cdot),\beta\right\rangle_t=\int_0^t\dfrac{d\left\langle x(u,\cdot),\xi\right\rangle_s}{\sigma(\xi_s)},\quad t\geq 0,\quad \text{a. s.}
$$
On the other hand,
$$
\left\langle x(u,\cdot),\xi\right\rangle_t=\left\langle x(u,\cdot),x(u,\cdot)-x(v,\cdot)\right\rangle_t=\left\langle x(u,\cdot)\right\rangle_t-\left\langle x(u,\cdot),x(v,\cdot)\right\rangle_t=
$$
$$
=t-\int_0^t\Phi(x(u,s)-x(v,s))ds=\dfrac{1}{2}\int_0^t\sigma^2(\xi_s)ds,\quad t\geq 0,\quad \text{a. s.}
$$
So,
$$
\left\langle\overline{x}(u,\cdot),\beta\right\rangle_t=\int_0^t\dfrac{\dfrac{1}{2}\sigma^2(\xi_s)}{\sigma(\xi_s)}ds=\dfrac{1}{2}\int_0^t\sigma(\xi_s)ds,\quad t\geq 0,\quad \text{a. s.}
$$
Thus, using Lemma 3.7, we get
$$
{\bf E} \left(\overline{x}(u,t)\vert\overline{x}(u,s)-\overline{x}(v,s), 0\leq s\leq t\right)=
$$
$$
={\bf E} \left(\overline{x}(u,t)\vert x(u,s)-x(v,s), 0\leq s\leq t\right)= {\bf E} \left(\overline{x}(u,t)\vert\xi_s, 0\leq s\leq t\right)=
$$
$$
=\int_0^t {\bf E} \left(\dfrac{d}{ds}\left\langle\overline{x}(u,\cdot),\beta\right\rangle_s|\xi_r, 0\leq r\leq s\right)d\beta_s=\dfrac{1}{2}\int_0^t\sigma(\xi_s)d\beta_s=
$$
$$
=\dfrac{1}{2}(\xi_t-(u-v))=\dfrac{1}{2}(\overline{x}(u,t)-\overline{x}(v,t)),\quad t\geq 0,\quad \text{a. s.},
$$
which proves $(3.6)$. Equality $(3.7)$ is a direct consequence of equality $(3.6)$. $\square$

\vskip10pt
\noindent {\bf Corollary 3.9.} {\it For any $u,v\in\mathbb{R}$ we have
$$
\lim_{t\rightarrow+\infty} {\bf E} \left[x(u,t)\Phi(x(u,t)-x(v,t))\right]= \dfrac{1}{2}(u+v).\eqno(3.8)
$$
}

\vskip10pt
\noindent {\bf Proof.} Using Theorem 3.4, we get the equalities
$$
{\bf E} \left[\overline{x}(u,t)\Phi(x(u,t)-x(v,t))\right]=
$$
$$
={\bf E} \left[\Phi(x(u,t)-x(v,t)) {\bf E} \left(\overline{x}(u,t)\vert\overline{x}(u,s)-\overline{x}(v,s), 0\leq s\leq t\right) \right]=
$$
$$
=\dfrac{1}{2} {\bf E} \left[\left(\overline{x}(u,t)-\overline{x}(v,t)\right)\Phi(x(u,t)-x(v,t))\right]=
$$
$$
=\dfrac{1}{2} {\bf E} \left[\left(x(u,t)-x(v,t)\right)\Phi(x(u,t)-x(v,t))\right]-
$$
$$
-\dfrac{1}{2}(u-v){\bf E} \Phi(x(u,t)-x(v,t)),
$$
and so,
$$
{\bf E} \left[x(u,t)\Phi(x(u,t)-x(v,t))\right]=
$$
$$
=\dfrac{1}{2} {\bf E} \left[\left(x(u,t)-x(v,t)\right)\Phi(x(u,t)-x(v,t))\right]+
$$
$$
+\dfrac{1}{2}(u+v){\bf E} \Phi(x(u,t)-x(v,t)).
$$
Applying $(2.2)$ and $(2.3)$ completes the proof. $\square$

\vskip10pt
\noindent {\bf Proposition 3.10.} {\it For any $u,v\in\mathbb{R}$ we have
$$
\lim_{t\rightarrow+\infty} \dfrac{1}{t} {\bf E} \left[x^2(u,t)x(v,t)\right]=u+2v.
$$
}

\vskip10pt
\noindent {\bf Proof.} Using It\^{o}'s formula and Fubini's theorem, we get
$$
{\bf E} \left[x^2(u,t)x(v,t)\right]=u^2v+vt+2\int_0^t {\bf E} \left[x(u,s)\Phi(x(u,s)-x(v,s))\right]ds,
$$
which together with $(3.8)$ yields the desired result. $\square$


\vskip30pt
\begin{thebibliography}{99}

\bibitem{Dean} Dean~D.~W. An analysis of the stochastic approaches to the problems of flow and transport in porous media// Thesis (PhD), University of Colorado at Denver, 1997. - 235~p. (electronic version)

\bibitem{Dorogovtsev2004} Dorogovtsev~A.~A. One Brownian stochastic flow// Theory of Stochastic Processes. - {\bf 10} ({\bf 26}), No.~3-4, 2004. - P.~21-25.

\bibitem{Dorogovtsev2007} Dorogovtsev~A.~A. Measure-valued processes and stochastic flows. - Kiev, the Institute of Mathematics of the NAS of Ukraine, 2007. - 289~p. (in Russian)

\bibitem{GrahamMcLaughlin1} Graham~W., McLaughlin~D. Stochastic analysis of nonstationary subsurface solute transport: 1.~Unconditional moments// Water Resources Research. - {\bf 25}, No.~2, 1989. - P.~215-232.

\bibitem{GrahamMcLaughlin2} Graham~W., McLaughlin~D. Stochastic analysis of nonstationary subsurface solute transport: 2.~Conditional moments// Water Resources Research. - {\bf 25}, No.~11, 1989. - P.~2331-2355.

\bibitem{HardyLittlewoodPolya} Hardy~G.~H., Littlewood~J.~E., P\'{o}lya~G. Inequalities. - Moscow, 1948. - 456~p. (in Russian)

\bibitem{Kotelenez} Kotelenez~P. A class of quasilinear stochastic partial differential equations of McKean-Vlasov type with mass conservation// Probability Theory and Related Fields. - {\bf 102}, No.~2, 1995. - P.~159-188.

\bibitem{Kunita} Kunita~H. Stochastic flows and stochastic differential equations. - Cambridge University Press, 1990. - 346~p.

\bibitem{LeJan} Le~Jan~Y. On isotropic Brownian motions// Zeitschrift f\"{u}r Wahrscheinlichkeitstheorie und verwandte Gebiete. - {\bf 70}, 1985. - P.~609-620.

\bibitem{LiptserShiryaev} Liptser~R.~S., Shiryaev~A.~N. Statistics of random processes (non-linear filtering and related topics). - Moscow, 1974. - 696~p. (in Russian)

\bibitem{Matsumoto} Matsumoto~H. Coalescing stochastic flows on the real line// Osaka Journal of Mathematics. - {\bf 26}, No.~1, 1989. - P.~139-158.

\bibitem{SotiropoulosKaznessis} Sotiropoulos~V., Kaznessis~Y. Analytical derivation of moment equations in stochastic chemical kinetics// Chemical Engineering Science. - {\bf 66}, No.~3, 2011. - P.~268-277.

\bibitem{Walsh} Walsh~J.~B. An introduction to stochastic partial differential equations. - In: \'{E}cole d'\'{e}t\'{e} de probabilit\'{e}s de Saint-Flour XIV--1984. Lecture Notes in Mathematics, {\bf 1180}. - Springer, Berlin, Heidelberg, 1986. - P.~265-439.

\bibitem{Zirbel} Zirbel~C.~L. Stochastic flows: dispersion of a mass distribution and Lagrangian observations of a random field// Dissertation (PhD), Princeton University, 1993. - 162~p. (electronic version)

\end{thebibliography}
\end{document}